\documentclass[12pt,leqno]{amsart}
\usepackage{amsmath,amsthm,amsfonts,amsbsy,multicol,fancybox}
\usepackage[utf8]{inputenc}
\usepackage[T1]{fontenc}

\usepackage{caption}\usepackage{subcaption}
\usepackage{graphicx}\usepackage{amssymb,latexsym}
\usepackage{epstopdf}
\usepackage{color}
\DeclareGraphicsExtensions{.pdf, .jpg, .pnp, .bmp, .fig}
\usepackage{tikz}
\usepackage[all]{xy}
\input xy
\xyoption{all}

\usepackage[hidelinks]{hyperref}
\usepackage{amssymb,latexsym}

\usepackage{thmtools}

\newcounter{argument}
\newenvironment{argument}[1][\medskip]{%
\refstepcounter{argument}
\par\medskip
\noindent\phantomsection
\textbf{#1~\thesection.\arabic{argument}\,\,}\rmfamily\em}{\hspace{\fill}$\Box$\par\smallskip\noindent}

\newcommand{\bass}{\begin{argument}[Assumption]}\newcommand{\eass}{\end{argument}}
\newcommand{\bth}{\begin{argument}[Theorem]} \newcommand{\ethe}{\end{argument}}
\newcommand{\bre}{\begin{argument}[Remark]}      \newcommand{\ere}{\end{argument}}
\newcommand{\ble}{\begin{argument}[Lemma]}       \newcommand{\ele}{\end{argument}}
\newcommand{\bde}{\begin{argument}[Definition]}   \newcommand{\ede}{\end{argument}}
\newcommand{\bco}{\begin{argument}[Corollary]}     \newcommand{\eco}{\end{argument}}
\newcommand{\bpr}{\begin{argument}[Proposition]}  \newcommand{\epr}{\end{argument}}
\newcommand{\bexam}{\begin{argument}[Example]}\newcommand{\eexam}{\end{argument}}

\newcommand{\bpf}{\begin{proof}}\newcommand{\epf}{\end{proof}}
\newcommand{\barr}{\begin{array}}\newcommand{\earr}{\end{array}}
\newcommand{\beao}{\begin{eqnarray*}}\newcommand{\eeao}{\end{eqnarray*}\noindent}
\newcommand{\beam}{\begin{eqnarray}}\newcommand{\eeam}{\end{eqnarray}\noindent}
\newcommand{\beqq}{\begin{equation}}\newcommand{\eeqq}{\end{equation}\noindent}

 \newcommand{\un}{\underbrace}

\newcommand{\nto}{n\to\infty}

\newcommand{\tto}{t\to\infty}

\newcommand{\D}{\Delta}

\newcommand{\w}{\omega}

\newcommand{\bfE}{{\mathbb E}}\newcommand{\bbE}{{\mathcal E}} 
\newcommand{\bbf}{{\mathcal F}}

\newcommand{\bbn}{{\mathcal N}} \newcommand{\bbN}{{\mathbb N}}

\newcommand{\bfP}{{\mathbb P}}

\graphicspath{{./../Images/}}

\begin{document}

\title[Pathwise stability and positivity of nonlinear SDEs]{A note on Pathwise stability and positivity of nonlinear stochastic differential equations}

\author[I. S. Stamatiou]{I. S. Stamatiou}
\email{joniou@gmail.com}

\begin{abstract}
We use the semi-discrete method, originally proposed in \emph{Halidias (2012), Semi-discrete approximations for stochastic differential equations and applications, International Journal of Computer Mathematics, 89(6)}, to reproduce qualitative properties of a class of nonlinear stochastic differential equations with nonnegative, non-globally Lipschitz coefficients and a unique equilibrium solution. The proposed fixed-time step method preserves the positivity of  solutions and reproduces the almost sure asymptotic stability behavior of the equilibrium with no time-step restrictions. 
\end{abstract}

\date\today

\keywords{Explicit Numerical Scheme; Semi-Discrete Method; non-linear SDEs; Stochastic Differential Equations; Boundary Preserving Numerical Algorithm; Pathwise Stability;
 \newline{\bf AMS subject classification 2010:}  60H10, 60H35, 65C20, 65C30, 65J15, 65L20.}
\maketitle

\section{Introduction}\label{PSP:sec:intro}
\setcounter{equation}{0}

We are interested in the following class of scalar stochastic differential equations (SDEs),
\beqq  \label{PSP-eq:scalarSDEs}
x_t =x_0 + \int_0^t x_sa(x_s)ds + \int_0^t x_sb(x_s)dW_s,
\eeqq
where $a(\cdot), b(\cdot)$ are non-negative functions with $b(u)\neq0$ for $u\neq0, x_0\geq0$ and
$\{W_{t}\}_{t\geq0}$ is a one-dimensional Wiener process adapted to the filtration  $\{\bbf_t\}_{t\geq0}.$ We want to reproduce dynamical properties of (\ref{PSP-eq:scalarSDEs}). We use a fixed-time step explicit numerical method, namely the semi-discrete method, which reads
\beqq  \label{PSP-eq:SDmethod}
y_{n+1} =y_n\exp\left\{\left(a(y_n)-\frac{b^2(y_n)}{2}\right)\D + b(y_n)\D W_n\right\}, \quad n\in\bbN,
\eeqq
with $y_0=x_0,$ where $\D=t_{n+1}-t_{n}$ is the time step-size and $\D W_n:= W_{t_{n+1}}- W_{t_{n}}$ are the increments of the Wiener process. For the derivation of (\ref{PSP-eq:SDmethod}) see Section \ref{PSP:sec:proofs}.

The scopes of this article are two. Our main goal is to reproduce the almost sure (a.s.) stability and instability of the unique equilibrium solution of (\ref{PSP-eq:scalarSDEs}), i.e. for the trivial solution $x_t\equiv0.$ The positivity of the drift pushes the solution to explosive situations and the diffusion stabilizes this effect in a way we want to mimic.

On the other hand, SDE (\ref{PSP-eq:scalarSDEs}) has unique positive solutions when $x_0>0.$ The semi-discrete method (\ref{PSP-eq:SDmethod}) preserves positivity by construction.

Explicit fixed-step Euler methods fail to strongly converge to solutions of (\ref{PSP-eq:scalarSDEs}) when the drift or diffusion coefficient grows superlinearly \cite[Theorem 1]{hutzenthaler_et_al.:2011}. Tamed Euler methods were proposed to overcome the aforementioned problem, cf. \cite[(4)]{hutzenthaler_jentzen:2015}, \cite[(3.1)]{tretyakov_zhang:2013}, \cite{sabanis:2016}
 and references therein; nevertheless in general they fail to preserve positivity. We also mention the method presented in \cite{neuenkirch_szpruch:2014} where they use the Lamperti-type transformation to remove the nonlinearity from the diffusion to the drift part of the SDE. Moreover, adaptive time-stepping strategies applied to explicit Euler method are an alternative way to address the problem and there is an ongoing research on that approach, see \cite{fand_giles:2016}, \cite{kelly_lord:2016} and  \cite{kelly_et_al:2017}. However, the fixed-step method we propose reproduces the almost sure asymptotic stability behavior of the equilibrium with no time-step restrictions, compare Theorems \ref{PSP-theorem:asstability} and  \ref{PSP-theorem:asinstability} with  \cite[Theorem 4.1 and 4.2]{kelly_et_al:2017} respectively. 
 
Our proposed fixed-step method is explicit, strongly convergent, non-explosive and positive. The semi-discrete method was originally proposed in \cite{halidias:2012} and further investigated  in  \cite{halidias_stamatiou:2016}, \cite{halidias:2014}, \cite{halidias:2015}, \cite{halidias:2015d}, \cite{halidias_stamatiou:2015} and \cite{stamatiou:2017}. We discretize in each subinterval (in an appropriate additive or multiplicative way) the drift and/or the diffusion coefficient producing a new SDE which we have to solve and not an algebraic equation as all the standard numerical methods. The way of discretization is implied by the form of the coefficients of the SDE and is not unique.
 
Let us now assume some minimal additional conditions for the functions $a(\cdot)$ and $b(\cdot).$ In particular, we assume locally Lipschitz continuity of  $a(\cdot)$ and $b(\cdot),$ which in turn implies the existence of a unique, continuous $\bbf_t$-measurable process $x$ (cf. \cite[Ch. 2]{mao:2007}) satisfying (\ref{PSP-eq:scalarSDEs}) up to the explosion time $\tau_e^{x_0},$ i.e. on the interval $[0,\tau_e^{x_0}),$ where 
$$
\tau_e^{x_0}:=\inf\{t>0: |x_t^{x_0}|\notin [0,\infty)\}.
$$
Denoting $\theta_e^{x_0}$ the first hitting time of zero, i.e. 
$$
\theta_e^{x_0}:=\inf\{t>0: |x_t^{x_0}|=0\},
$$
it was shown in \cite[Section 3]{appleby_et_al:2008} that in the case
\beqq  \label{PSP-eq:condition}
\sup_{u\neq0} \frac{2a(u)}{b^2(u)}=\beta<1,
\eeqq
then $\tau_e^{x_0}=\theta_e^{x_0}=\infty,$ i.e. there exist unique positive solutions. The equilibrium zero solution of (\ref{PSP-eq:scalarSDEs}) is a.s. stable if  (see again \cite[Section 3]{appleby_et_al:2008})
\beqq  \label{PSP-eq:condition2}
\lim_{u\rightarrow0} \frac{2a(u)}{b^2(u)}<1,
\eeqq
i.e. for all $x_0>0$ 
$$
\bfP (\{\w: \lim_{\tto}x_t(\w)=0\})>0.
$$
Condition (\ref{PSP-eq:condition2}) shows how condition (\ref{PSP-eq:condition}) is close to being sharp. Furthermore, the presence of a sufficiently intense stochastic perturbation (because of the positivity of the function $a(\cdot)$) is necessary for the existence of a unique global solution and stability of the zero equilibrium.

The outline of the article is the following. In Section \ref{PSP:sec:main} we present our main results, that is Theorems \ref{PSP-theorem:asstability} and \ref{PSP-theorem:asinstability}, the proofs of which are deferred to Section \ref{PSP:sec:proofs}. Section \ref{PSP:sec:numerics} provides a numerical example.

\section{Main results}\label{PSP:sec:main}
\setcounter{equation}{0}

\bass\label{NSF:assA} 
The functions  $a(\cdot)$ and $b(\cdot)$ of (\ref{PSP-eq:scalarSDEs})  are non-negative with $b(u)\neq0$ for $u\neq0.$ 
\eass

The following result provides sufficient conditions for solutions of the semi-discrete scheme (\ref{PSP-eq:SDmethod}) to demonstrate a.s. stability.

\bth\label{PSP-theorem:asstability}[a.s. stability]
Let $a(\cdot)$ and $b(\cdot)$ satisfy Assumption \ref{NSF:assA} and (\ref{PSP-eq:condition}), i.e. there exists $\beta<1$ such that
\beqq  \label{PSP-eq:conditionsup}
\sup_{u\neq0} \frac{2a(u)}{b^2(u)}=\beta.
\eeqq
Let also $\{y_n\}_{n\in\bbN}$ be a solution of (\ref{PSP-eq:SDmethod})  with $y_0=x_0>0.$
Then for all $\D>0$,
\beqq  \label{PSP-eq:asstability}
\lim_{\nto} y_n=0, \qquad \hbox{a.s.}
\eeqq 
\ethe

The following result provides sufficient conditions for solutions of the semi-discrete scheme (\ref{PSP-eq:SDmethod}) to demonstrate a.s. instability. 

\bth\label{PSP-theorem:asinstability}[a.s. instability]
Let $a(\cdot)$ and $b(\cdot)$ satisfy Assumption \ref{NSF:assA} and there exists $\gamma>1$ such that
\beqq  \label{PSP-eq:conditioninf}
\liminf_{u\rightarrow0} \frac{2a(u)}{b^2(u)}=\gamma.
\eeqq
Let also $\{y_n\}_{n\in\bbN}$ be a solution of (\ref{PSP-eq:SDmethod})  with $y_0=x_0>0.$
Then for all $\D>0$,
\beqq  \label{PSP-eq:asinstability}
\bfP (\{\w: \lim_{\nto} y_n(\w)=0\})=0.
\eeqq 
\ethe

Note that there is no time-step restriction in the results of Theorems \ref{PSP-theorem:asstability} and \ref{PSP-theorem:asinstability}, i.e. (\ref{PSP-eq:asstability}) and (\ref{PSP-eq:asinstability})  hold for all $\D>0.$

\section{Numerical illustration}\label{PSP:sec:numerics}
\setcounter{equation}{0}

\begin{figure}[ht]
\centering
\begin{subfigure}{.6\textwidth}
   \includegraphics[width=1\textwidth]{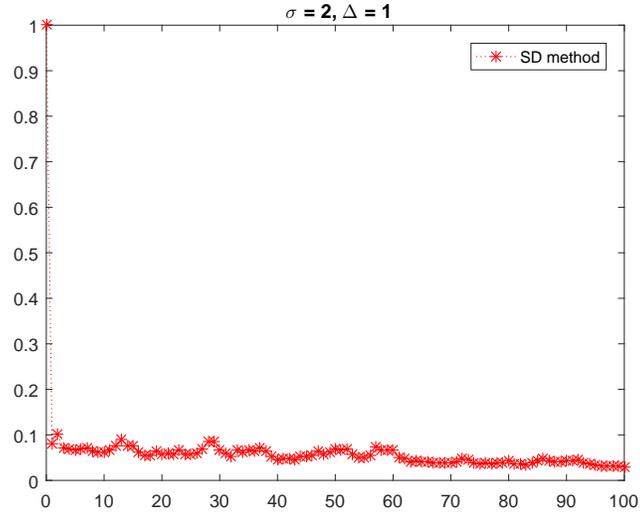}\label{PSP-fig:s2}
     \caption{Trajectory for (\ref{PSP-eq:exampleSD}) with  $\sigma=2$.}
  \end{subfigure}
  \begin{subfigure}{.6\textwidth}
  \centering
   \includegraphics[width=1\textwidth]{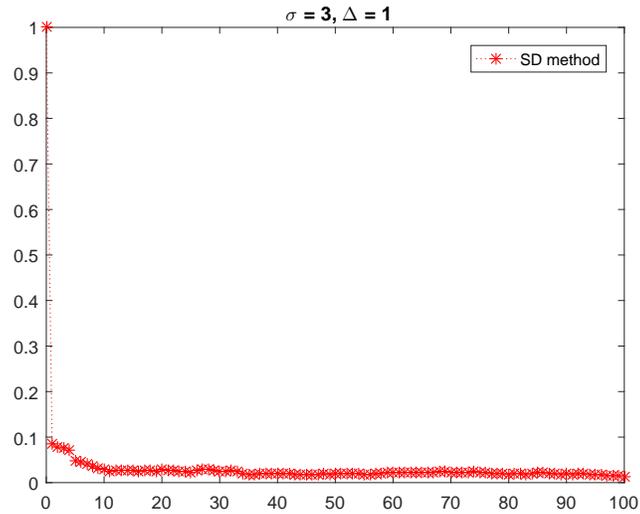}\label{PSP-fig:s3}
     \caption{Trajectory for (\ref{PSP-eq:exampleSD}) with  $\sigma=3.$}
\end{subfigure}
 \caption{Trajectories of (\ref{PSP-eq:exampleSD}) for different values of $\sigma$.}\label{PSP-fig:s23}
  \end{figure}

We will use the numerical example of \cite[Section 5]{kelly_et_al:2017}, that is we take $a(x)=x^2$ and $b(x)=\sigma x$ and $x_0=1$ in (\ref{PSP-eq:scalarSDEs}), i.e.
\beqq  \label{PSP-eq:exampleSDE}
x_t =1 + \int_0^t (x_s)^3ds + \sigma\int_0^t (x_s)^2dW_s, \qquad t\geq0.
\eeqq
Note that the value of $\sigma$ determines the value of the ratio $2a(u)/b^2(u)=2/\sigma^2.$ The semi-discrete method (\ref{PSP-eq:SDmethod}) reads
\beqq  \label{PSP-eq:exampleSD}
y_{n+1} =y_n\exp\left\{\left(1-\frac{\sigma^2}{2}\right)(y_n)^2\D + \sigma y_n\D W_n\right\}, \quad n\in\bbN,
\eeqq
with $y_0=1.$ First we examine the case of stability, that is when $\beta<1$ or $\sigma>\sqrt{2}.$ Figure \ref{PSP-fig:s23} displays a trajectory of the semi-discrete method (\ref{PSP-eq:exampleSD}) for the cases $\sigma=2$ and $\sigma=3$ accordingly. We observe the asymptotic stability in each case as well as the positivity of the paths. There is no need for time step restriction as in  \cite[Fig. 2]{kelly_et_al:2017}.

\begin{figure}[ht]
\centering
\begin{subfigure}{.6\textwidth}
   \includegraphics[width=1\textwidth]{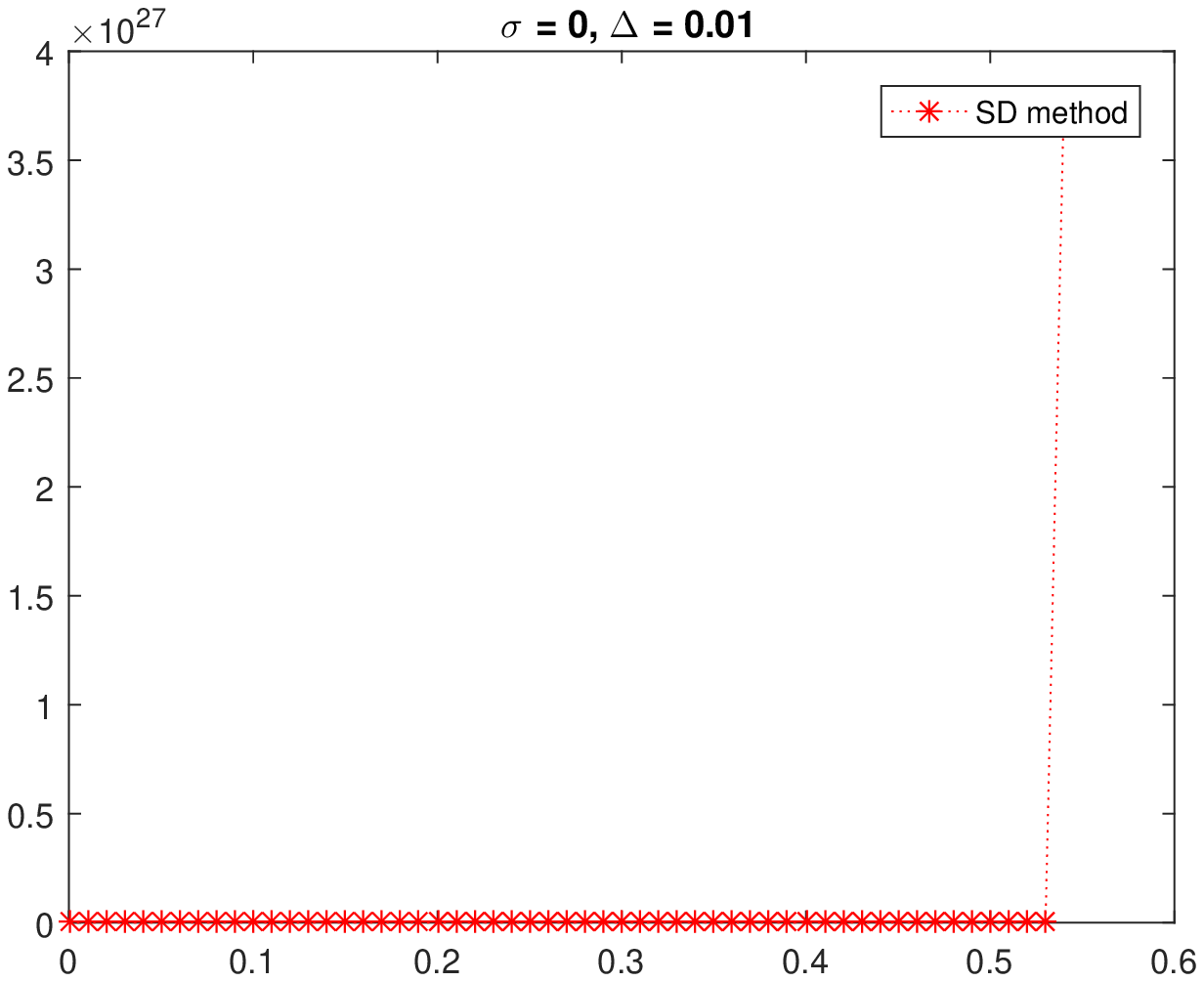}\label{PSP-fig:s0}
     \caption{Trajectory for (\ref{PSP-eq:exampleSD}) with  $\sigma=0$ and $\D=0.01.$}
  \end{subfigure}
  \begin{subfigure}{.6\textwidth}
  \centering
   \includegraphics[width=1\textwidth]{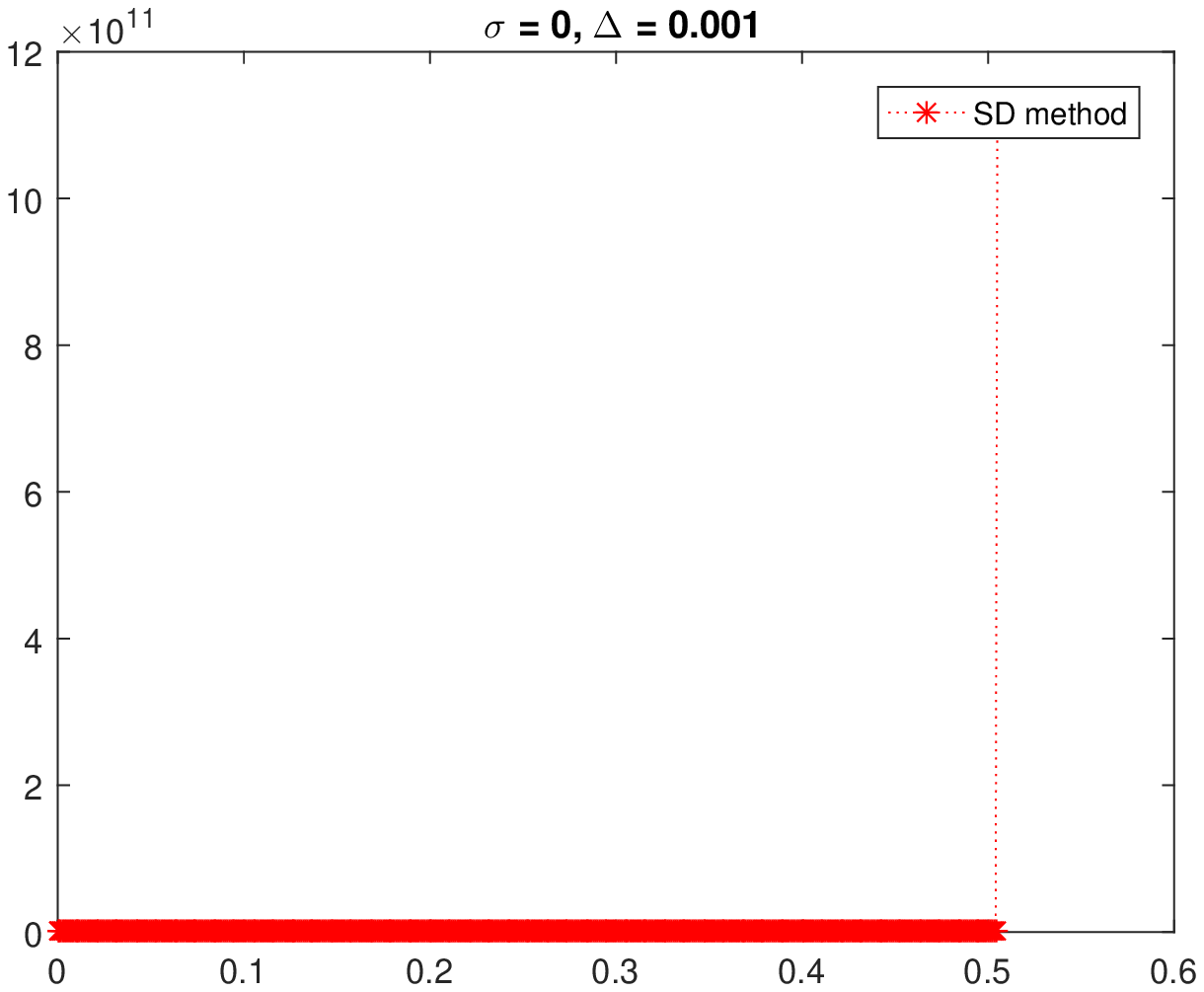}\label{PSP-fig:s0b}
     \caption{Trajectory for (\ref{PSP-eq:exampleSD}) with  $\sigma=0$ and $\D=0.001.$}
\end{subfigure}
 \caption{Trajectories of (\ref{PSP-eq:exampleSD}) for $\sigma=0$ and different values of $\D.$.}\label{PSP-fig:s00b}
  \end{figure}

\begin{figure}[ht]
\centering
\begin{subfigure}{.6\textwidth}
   \includegraphics[width=1\textwidth]{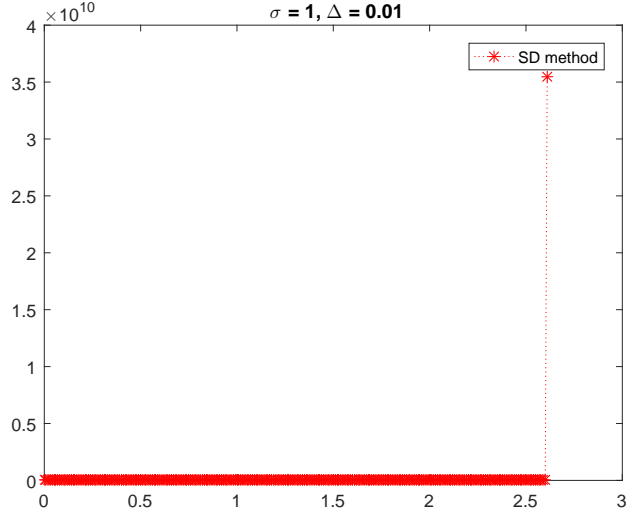}\label{PSP-fig:s1}
     \caption{Trajectory for (\ref{PSP-eq:exampleSD}) with  $\sigma=1$ and $\D=0.01.$}
  \end{subfigure}
  \begin{subfigure}{.6\textwidth}
  \centering
   \includegraphics[width=1\textwidth]{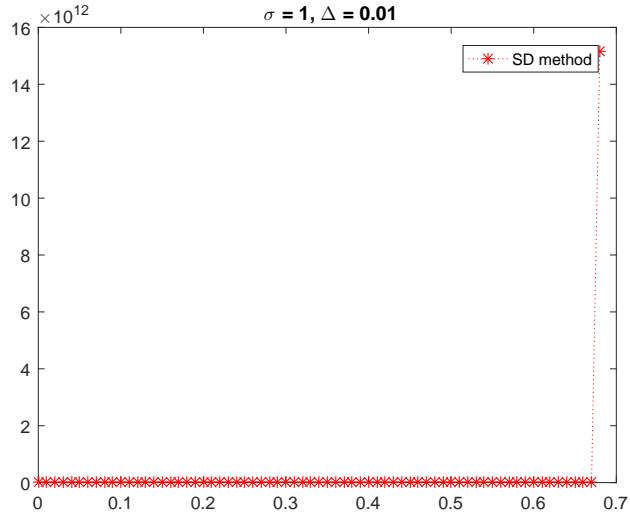}\label{PSP-fig:s1b}
     \caption{Trajectory for (\ref{PSP-eq:exampleSD}) with  $\sigma=1$ and $\D=0.01.$}
\end{subfigure}
 \caption{Trajectories of (\ref{PSP-eq:exampleSD}) for $\sigma=1.$}\label{PSP-fig:s11b}
  \end{figure}

 Figures \ref{PSP-fig:s00b} and  \ref{PSP-fig:s11b} displays the case when $\gamma>1$ or  $\sigma<\sqrt{2}.$ We consider the cases  $\sigma=0$ and $\sigma=1$ accordingly. Now, we observe instability and an apparent finite-time explosion. The apparent explosion time in the ordinary differential equation (case $\sigma=0$)  is very close to the computed one
 $$
 \tau_e^1:= \int_1^\infty u^{-3}du=0.5,
 $$
 and becomes closer as we lower the step-size $\D.$ In the case $\sigma=1$ we observe again  the apparent explosion time for the SDE which is now random.

\section{Proofs}\label{PSP:sec:proofs}
\setcounter{equation}{0}

In this section we first discuss about the derivation of the semi-discrete scheme (\ref{PSP-eq:SDmethod}) and then we provide the proofs of Theorems \ref{PSP-theorem:asstability} and \ref{PSP-theorem:asinstability}.

Given the equidistant partition $0=t_0<t_1<\ldots<t_N=T$ with step size $\D=T/N$ we consider the following process
\beqq  \label{PSP-eq:SDprocess}
y_{t} =y_{t_n} + \int_{t_n}^t \un{a(y_{t_n})y_s}_{f(y_{t_n},y_s)}ds + \int_{t_n}^t \un{b(y_{t_n})y_s}_{g(y_{t_n},y_s)} dW_s,
\eeqq
in each subinterval $[t_n,t_{n+1}],$ with $y_0=x_0$ a.s. By (\ref{PSP-eq:SDprocess}) the form of discretization becomes apparent. We discretized the drift and diffusion coefficient in a multiplicative way producing a new SDE at each subinterval with the unique strong solution
\beqq  \label{PSP-eq:SDsolutionprocess}
y_{t} =y_{t_n}\exp\left\{\left(a(y_{t_n})-\frac{b^2(y_{t_n})}{2}\right)(t-t_n) + b(y_{t_n})(W_t- W_{t_n})\right\}.
\eeqq
The first variable of the auxiliary functions $f(\cdot,\cdot)$ and $g(\cdot,\cdot)$ in (\ref{PSP-eq:SDprocess}) denote the discretized part. In case $a(\cdot)$ and $b(\cdot)$ are locally Lipschitz so are $f(\cdot,\cdot)$ and $g(\cdot,\cdot)$ and as a consequence we have a strong convergence result of the type (see \cite[Theorem 2.1]{halidias_stamatiou:2016})
$$
\lim_{\D\rightarrow0}\bfE \sup_{0\leq t\leq T}|y_t-x_t|^2=0,
$$
in the case of finite moments of the original SDE and the approximation process, that is when $\bfE|x_t|^p\vee\bfE|y_t|^p<A$ for some $p>2$ and $A>0.$ However, the strong convergence of the method  does not hold in all cases considered here since the moments of $(x_t)$ are bounded only up to an explosion time $\tau_e^{x_0}$ which may be finite in case (\ref{PSP-eq:condition}) does not hold. The main focus here is the preservation of the dynamics in the discretization as shown in Section \ref{PSP:sec:numerics}.

The stability behavior of the equilibrium solution of (\ref{PSP-eq:SDmethod}) is an easy task since we have an analytic expression of the solution process. Nevertheless, we discuss the steps below.

Take a $p>0$ to be specified later on and rewrite (\ref{PSP-eq:SDmethod})  as
\beao
|y_{n+1}|^p &=&|y_{n}|^p\exp\left\{\frac{pb^2(y_n)}{2}\left(\frac{2a(y_n)}{b^2(y_{n})}-1+p\right)\D\right\}\exp\left\{-\frac{p^2b^2(y_{n})}{2}\D + pb(y_{n})\D W_{n}\right\}\\
&=&\bbE(y_{n})\xi_{n+1},
\eeao
where we used the notation $y_n$ for  $y_{t_n},$ the exponential function $\bbE(\cdot)$ reads
$$
\bbE(u) = \exp\left\{\frac{pb^2(u)}{2}\left(\frac{2a(u)}{b^2(u)}-1+p\right)\D\right\}
$$
and for $t\in (t_n,t_{n+1}]$ we consider the SDE
$$
d\xi_t = pb(y_n)dW_t
$$
with $\xi_n$=$|y_{n}|^p.$ Therefore $\bfE\xi_{n+1}$=$\bfE|y_{n}|^p$ and choosing $0<p<1-\beta,$ where $\beta$ is as in Theorem \ref{PSP-theorem:asstability}, we get that $\bbE(u)\leq1$ for any $\D>0$ implying on the one hand the boundness of the moments of ${y_n}_{n\in\bbn}$ and on the other hand the $p_th$ moment exponential stability of the trivial solution of $y_t^{x_0}.$ This in turn implies the a.s. exponential stability of the trivial solution (see \cite[Theorem 4.4.2]{mao:2007}) and consequently (\ref{PSP-eq:asstability}). As a result we also now the rate of (\ref{PSP-eq:asstability}) which is exponential and determined by the function $\bbE(\cdot).$ The result of Theorem \ref{PSP-theorem:asinstability} follows by analogue arguments where now we consider the representation
\beao
|y_{n+1}|^{-p} &=&|y_{n}|^{-p}\exp\left\{\frac{pb^2(y_n)}{2}\left(\frac{-2a(y_n)}{b^2(y_{n})}+1+p\right)\D\right\}\exp\left\{-\frac{p^2b^2(y_{n})}{2}\D - pb(y_{n})\D W_{n}\right\}\\
&=&\bbE^*(y_{n})\xi^*_{n+1},
\eeao
for $0<p<\gamma-1$ where $\gamma$ is as in the statement of Theorem \ref{PSP-theorem:asinstability} and 
$$
\bbE^*(u) = \exp\left\{\frac{pb^2(u)}{2}\left(-\frac{2a(u)}{b^2(u)}+1+p\right)\D\right\}
$$
and for $t\in (t_n,t_{n+1}]$ 
$$
d\xi^*_t = -pb(y_n)dW_t.
$$

 

\bibliographystyle{unsrt}\baselineskip12pt 
\bibliography{stability_positivity}

\appendix

\end{document}